\documentclass[preprint,3p,12pt]{elsarticle}

\usepackage{lineno,hyperref}
\modulolinenumbers[5]

\usepackage{hyperref}

\usepackage{epstopdf}

\usepackage{mathrsfs}
\usepackage{amssymb}
\usepackage{amsfonts}
\usepackage{latexsym}
\usepackage{amsmath}
\usepackage{xcolor}
\usepackage{wrapfig}
\usepackage{floatflt}

\usepackage{mathtools}
\usepackage{extarrows}

\usepackage{graphicx}
\def\lb{\label}

\newcommand{\er}[1]{\textrm{(\ref{#1})}}

\newtheorem{theorem}{\bf Theorem}[section]

\def\a{\alpha}         
\def\b{\beta}

\def\d{\delta}

\def\r{\rho}           
\def\s{\sigma}         
\def\S{\Sigma}         
\def\t{\tau}

\def\O{\Omega}

\def\ve{\varepsilon}   \def\vt{\vartheta}        

       \def\C{{\mathbb C}}    
    \def\N{{\mathbb N}}   

\def\lt{\biggl}                  \def\rt{\biggr}
               \def\wt{\widetilde}
\def\no{\noindent}

\let\ge\geqslant                 \let\le\leqslant

\def\iy{\infty}
\def\sm{\setminus}               \def\es{\emptyset}
\def\ss{\subset}                 \def\ts{\times}

\def\el2{\ell^{\,2}}             \def\1{1\!\!1}

\def\Im{\mathop{\mathrm{Im}}\nolimits}

\def\Re{\mathop{\mathrm{Re}}\nolimits}

\def\supp{\mathop{\mathrm{supp}}\nolimits}

\newtheorem{corollary}[theorem]{\bf Corollary}

\let\ge\geqslant
\let\le\leqslant

\newcommand{\ca}{\begin{cases}}
\newcommand{\ac}{\end{cases}}
\newcommand{\ma}{\begin{pmatrix}}
\newcommand{\am}{\end{pmatrix}}
\def\eq{\begin{equation}}
\def\qe{\end{equation}}
\def\[{\begin{equation}}
\def\]{\end{equation}}

\begin{document}

\begin{frontmatter}




\title{Vi\`ete's fractal distributions and their momenta}

\date{\today}

\author
{Anton A. Kutsenko}

\address{Jacobs University, 28759 Bremen, Germany; email: akucenko@gmail.com}

\begin{abstract}
Solutions of Schr\"oder-Poincar\'e's polynomial equations $f(az)=P(f(z))$ usually do not admit a simple closed-form representation in terms of known standard functions. We show that there is a one-to-one correspondence between zeros of $f$ and a set of discrete functions stable at infinity. The corresponding Vi\`ete-type infinite products for zeros of $f$ are also provided. This allows us to obtain a special kind of closed-form representation for $f$ based on the Weierstrass-Hadamard factorization. From this representation, it is possible to derive explicit momenta formulas for zeros. We discuss also the rate of convergence of WH-factorization and momenta formulas. Obtaining explicit closed-form expressions is the main motivation for this work. Finally, all the branches of the multi-valued function $f^{-1}$ are computed explicitly.  
\end{abstract}




\end{frontmatter}

{\section{Introduction and main results}\lb{sec0}}

The classical Vi\`ete's formula
$$
 \frac{2}{\pi}=\frac{\sqrt{2}}{2}\cdot\frac{\sqrt{2+\sqrt{2}}}{2}\cdot\frac{\sqrt{2+\sqrt{2+\sqrt{2}}}}{2}...
$$
uses nested square root radicals to represent the constant $\pi$. Wiki says "By now many formulas similar to Vi\`ete's involving either nested radicals or infinite products of trigonometric functions are known for $\pi$, as well as for other constants such as the golden ratio", see, e.g., \cite{S,N,P,L}. In this note, we derive formulas for zeros of functions satisfying Schr\"oder-Poincar\'e's polynomial equations. In general, the formulas for zeros will involve various nested-radicals products similar to Vi\`ete's. These formulas can be used in Weierstrass-Hadamard factorization to obtain various closed-form expressions.

Finally, looking through "A chronology of continued square roots and other continued compositions" \cite{J}, I found paper \cite{AK}, where a detailed analysis of real roots of $f$, satisfying $f(az)=f(z)^2+c$, is provided. Many interesting facts are presented in \cite{J}, e.g., an interesting story of the famous formula
$$
 \ve_0\sqrt{2+\ve_1\sqrt{2+\ve_2\sqrt{2+...}}}=2\sin\lt(\frac{\pi}{4}\sum_{n=0}^{+\iy}\frac{\ve_0\ve_1\ve_2...\ve_n}{2^n}\rt),
$$
where $\ve_i=-1,0,1$.

We assume facts about existence of entire solutions of SP-equation to be known, see, e.g., \cite{F1,DGV}. Let $P$ be some polynomial of degree $d\ge2$. Let $b$ be some its repelling point $P(b)=b$, with $|a|>1$ for $a:=P'(b)$. Consider the entire solution $f$ of SP-equation $f(az)=P(f(z))$ satisfying $f(0)=b$, $f'(0)=1$. 
This solution can be taken as
\[\lb{limit}
 f(z)=\lim_{n\to\iy}\underbrace{P\circ...\circ P}_n(b+a^{-n}z),\ \ z\in\C,
\] 
see, e.g., \cite{G}. Composition \er{limit} converges uniformly in any compact subset of $\C$. For simplicity, let us assume $b\neq0$. This is not a restriction, since $\wt f:=f+c$, $c\in\C$, also satisfies some polynomial SP-equation. Let $P_0^{-1}(w)$, $w\in\C$ be the principal branch of $P^{-1}$ analytic in some open domain containing $b$, where $P_0^{-1}(b)=b$. We assume also that  

{\no\bf Hypothesis 1.} For any $w\in\C$ the orbit 
\[\lb{H1}
 (P_0^{-1})^{\circ n}(w):=\underbrace{P_0^{-1}\circ...\circ P_{0}^{-1}}_n(w)\to b.
\]

This assumption means that point $b$ repelling for $P$ is attracting for $P_0^{-1}$. Note that once $(P_0^{-1})^{\circ k}(w)\in \{|w-b|<\d\}$ for some small $\d> 0$ and some $k\in\N$, then $(P_0^{-1})^{\circ n}(w)$ stays in $\{|w-b|<\d\}$ for $n>k$ and $(P_0^{-1})^{\circ n}(w)\to b$, since 
$$
 |(P_0^{-1})'(b)|=|a^{-1}|<1.
$$
Let $P_j^{-1}$, $1\le j\le d-1$ be other branches of $P^{-1}$ so that $\{z_j(w)\}_{j=0}^{d-1}=\{P_j^{-1}(w)\}_{j=0}^{d-1}$ is a complete set of solutions of $P(z)=w$, defined for all $w\in\C$. For our research, it does not matter how the branches of $P^{-1}$ are numbered. There are only two things that we should pay close attention to: 1) analyticity of the principal branch $P_0^{-1}$ at an open neighbourhood of its attracting point $b$; 2) {\bf Hypothesis 1}.

Introduce the polynomial
\[\lb{Q}
 Q(z):=\frac{P(z)-P(b)}{z-b}=\frac{P(z)-b}{z-b}
\]
and the set of discrete functions stable at infinity
\[\lb{S}
 \S=\{\s:\N\to\{0,...,d-1\},\ \lim_{n\to\iy}\s_n=0\}.
\] 

\begin{theorem}\lb{T1} The set of zeros of $f$ coincides with $\{z(\s)\}_{\s\in\S}$, where
\[\lb{all_zeros}
z(\s)=-b\prod_{n=1}^{\iy}\frac{a}{Q(P_{\s_n}^{-1}\circ...\circ P_{\s_1}^{-1}(0))}.
\]
Each zero is counted according to its multiplicity. In other words, the multiplicity of $z_0$ as zero of $f$ is $\#\{\s\in\S:\ z(\s)=z_0\}$.
\end{theorem}

We may apply Theorem \ref{T1} to the function $F(z):=f(z)-w$ with some constant $w\in\C$, because $F$ also satisfies SP-equation $F(az)=P(F(z)+w)-w$ similar to that for $f$. We only should care about the assumption $f(0)-w=b-w\ne0$, see after \er{limit}.

\begin{corollary}\lb{C1}
All the solutions of $f(z)=w$, where $w\in\C\sm\{b\}$, have the form
\[\lb{branchesg}
 z=g_{\s}(w)=(w-b)\prod_{n=1}^{\iy}\frac{a}{Q(P^{-1}_{\s_n}\circ...\circ P^{-1}_{\s_1}(w))},\ \ \s\in\S.
\]
Each solution is counted according to its multiplicity.
\end{corollary}

The case $w=b$ is special, because formal setting $w=b$ in \er{branchesg} leads to $g_{\s}(b)=0$, $\forall\s\in\S$ which seems doubt. The next theorem is devoted to the case $w=b$.

\begin{theorem}\lb{T2}
All the solutions of $f(z)=b$ have the form $z=g_{0,0}=0$ or
\[\lb{special}
 z=g_{\s,m}=a^m(P_{\s_1}^{-1}(b)-b)\prod_{n=2}^{\iy}\frac{a}{Q(P^{-1}_{\s_n}\circ...\circ P^{-1}_{\s_1}(b))},\ \ m\in\N\cup\{0\},\ \ \s\in\S',
\]
where $\S'=\{\s\in\S:\ \s_1\ne0\}$.
\end{theorem}

In fact, $\{g_{\s}\}_{\s\in\S}$ are all the branches of super-multi-valued function $f^{-1}$. Depending on the choice of the branches $P^{-1}$, the functions $g_{\s}(w)$ may or may not be analytic. We can only state that the branch $g_0(w)$ ($\s=0$) is analytic in $\O\sm\{b\}$, where $\O$ is some small neighborhood of $b$, e.g., considered in the remark after {\bf Hypothesis 1}. Due to \er{limit} and arguments presented before \er{e4}, we have
\[\lb{g0}
 g_0(w)=\lim_{n\to\iy}a^n(\underbrace{P_0^{-1}\circ...\circ P_0^{-1}}_n(w)-b).
\]     
Differentiating \er{g0} and using $(P_0^{-1})'(w)=1/P'(P_0^{-1}(w))$, $w\in\O$ we obtain
\[\lb{g0p}
 g_0'(w)=\prod_{n=1}^{\iy}\frac{a}{P'(\underbrace{P_0^{-1}\circ...\circ P_0^{-1}}_n(w))}.
\]
The convergence of the product \er{g0p}, as well as \er{g0}, and \er{branchesg}, \er{all_zeros} is exponentially fast, as discussed in the beginning of Section \ref{Proof}. 

The order of the entire function $f$ can be computed explicitly by substituting $\a e^{A|z|^{\r}}$ into SP-equation $f(az)=P(f(z))$, see, e.g., \cite{M}. Extracting leading terms after the substitution, we obtain $\r=\frac{\ln d}{\ln |a|}$. If $\r<1$ then the Weierstrass-Hadamard (WH) factorization for $f$ does not contain exponential factors, see \cite{T}.

\begin{corollary}\lb{C2}
Suppose that $d<|a|$. If $w\ne b$ then WH-factorization for $f$ is
\[\lb{gen}
 f(z)=w+(b-w)\prod_{\s\in\S}\lt(1-\frac{z}{g_{\s}(w)}\rt),\ \ z\in\C.
\]
In particular,
\[\lb{WHgen}
 f(z)=b\prod_{\s\in\S}\lt(1+\frac{z}{b}\prod_{n=1}^{\iy}\frac{Q(P_{\s_n}^{-1}\circ...\circ P_{\s_1}^{-1}(0))}{a}\rt),\ \ z\in\C.
\] 
If $w=b$ then
\[\lb{WHspecial}
 f(z)=b+z\prod_{m\ge0}\prod_{\s\in\S'}\lt(1-\frac{z}{a^m g_{\s,0}}\rt)=b+\lt(\lt\{\frac{z}{g_{\s,0}}\rt\}_{\s\in\S'};\frac1a\rt)_{\iy}z,\ \ z\in\C,
\]
where $(\{\a\}_{i\in R};\b)_{\iy}:=\prod_{i\in R}\prod_{n=0}^{\iy}(1-\a_i\b^n)$ is the q-Pochhammer symbol.
\end{corollary}

Equation \er{gen} allows us to compute explicitly momentum formulas for zeros of $f(z)-w$, for any fixed $w\ne b$. This new type of formulas will include both: infinite products and infinite sums. The first (negative) momentum formula for zeros follows from \er{gen} immediately 
\[\lb{moment}
 \sum_{\s\in\S}\prod_{n=1}^{\iy}\frac{Q(P_{\s_n}^{-1}\circ...\circ P_{\s_1}^{-1}(w))}{a}=f'(0)=1,\ \ \forall w\in\C\sm\{b\}.
\]

Let us note how to compute explicitly other momenta of zeros. First, differentiating $f(az)=P(f(z))$ at $z=0$ and using $f(0)=b$, $f'(0)=1$, $P'(b)=a$, we obtain recurrent formulas to determine all the derivatives: 
\begin{align} 
 f''(0)&=(a^2-a)^{-1}P''(b),\lb{d2}  \\ 
 f^{(m)}(0)&=(a^m-a)^{-1}\sum_{j=2}^{m}P^{(j)}(b)B_{m,j}(f'(0),...,f^{(m-j+1)}(0)),\ \ m\ge2,\lb{dm}
\end{align}
where $B_{m,j}$ are Bell polynomials. They are given by
\[\lb{Bell}
 B_{m,j}(x_1,...,x_{m-j+1})=\sum\frac{m!}{k_1!...k_{m-j+1}!}\lt(\frac{x_1}{1!}\rt)^{k_1}...\lt(\frac{x_{m-j+1}}{(m-j+1)!}\rt)^{k_{m-j+1}},
\] 
where the sum is taken over all sequences $k_1$, $k_2$, ..., $k_{m-j+1}$ of non-negative integers such that the two conditions are satisfied:
\[\lb{Bell1}
 \sum_{i=1}^{m-j+1}k_i=j,\ \ \sum_{i=1}^{m-j+1}ik_i=m,
\]
see more about Fa\`a di Bruno's formula for high order derivatives of compositions in, e.g., wiki. Now, differentiating $\ln (f(z)-w)$ at $z=0$ and using \er{WHgen}, we obtain the momenta formulas of high orders $m\ge2$:
\begin{align} 
\sum_{\s\in\S}\prod_{n=1}^{\iy}\frac{Q(P_{\s_n}^{-1}\circ...\circ P_{\s_1}^{-1}(w))^2}{a^2}&=f'(0)^2-bf''(0)=1-\frac{b P''(b)}{a^2-a}\lb{moment2} \\
\sum_{\s\in\S}\prod_{n=1}^{\iy}\frac{Q(P_{\s_n}^{-1}\circ...\circ P_{\s_1}^{-1}(w))^m}{a^m}&=\sum_{j=1}^{m}\frac{(-b)^{m-j}(j-1)!}{(m-1)!}B_{m,j}(f'(0),...,f^{(m-j+1)}(0)).\lb{momentm}
\end{align}
Another type of  
Vieta formulas also follows from \er{gen}:
\[\lb{Vgen}
 \sum_{\s\in\S}\frac{w-b}{g_{\s}(w)}=f'(0)=1,\ \ \ 
 \sum_{\s\ne\t}\frac{(w-b)^2}{g_{\s}(w)g_{\t}(w)}=\frac{f''(0)}{2!}=\frac{P''(b)}{2(a^2-a)}
\]
and so on. There is a natural extension of $g_{\s}(w)$ to the case $w=b$:
\[\lb{gsgen}
 g_{\s}(b)=\ca 0, & \s=0,\\
               g_{\s',m}, & \s'=\s|_{\{n:\ n\ge m\}},\ {\rm where}\ \s_m\ne0\ {\rm and}\ \s_n=0\ {\rm for}\ n<m,\ac\ \ \s\in\S.
\] 
Then {\it all the solutions $f(z)=w$, where $w\in\C$, have the form $z=g_{\s}(w)$, $\s\in\S$}. It is possible to estimate the remainder of series \er{momentm}, \er{moment}, and products \er{gen}-\er{WHspecial}. They converge exponentially  fast, regarding the length of the support of functions $\s\in\S$.

\begin{theorem}\lb{T3}
There is $C>0$, depending on the polynomial $P$ only, such that
\[\lb{est1}
 |g_{\s}(w)|\ge C|a|^{|\supp\s|},\ \s\in\S,\ \ \s\ne0.
\]
Moreover, if $d<|a|$ then
\[\lb{est2}
 \sum_{|\supp\s|\ge N}|g_{\s}(w)|^{-m}=O((d|a|^{-m})^N),\ \ \prod_{|\supp\s|\ge N}\lt|1-\frac{z}{g_{\s}(w)}\rt|=1+O(|z|(d|a|^{-1})^N),
\]
for $N\to\iy$, $m\ge1$, and any fixed $z\in\C$. 
\end{theorem}

We will begin the next section with examples. The proof of the main result is placed in the final section.

{\section{Examples}\lb{sec1}} 

{\bf 1.} Consider the case $P(z)=2z^2-1$. SP-equation is $f(az)=2f(z)^2-1$. We take $f(0)=b=1$, $f'(0)=1$. Then $a=(2z^2-1)'|_{z=b}=4$. Polynomial \er{Q} is 
$$
 Q(z)=\frac{2z^2-1-1}{z-1}=2z+2.
$$
There are two branches of $P^{-1}$:
$$
 P_{1}^{-1}(w)=\sqrt{\frac{1+w}2},\ \ P_{-1}^{-1}(w)=-\sqrt{\frac{1+w}2}.
$$
We assume that
$$
 \sqrt{z}=r^{\frac12}e^{\frac{i\vt}2}\ \ {\rm for}\ \ z=r e^{i\vt},\  r\ge0,\ \vt\in(-\pi,\pi].
$$
The branch $P_1^{-1}$ is principal. It is analytically defined near the attracting (for $P_1^{-1}$) point $b$. Moreover, $(P_1^{-1})^{\circ n}(w)$ converges to its fixed point $b$ for any $w\in\C$, since $\sqrt{z}$ is a contraction mapping in the closed domain $D=\{z:\ \Re z\ge1/\sqrt{2}\}$:
$$
 |\sqrt{z_1}-\sqrt{z_2}|=\frac{|z_1-z_2|}{|\sqrt{z_1}+\sqrt{z_2}|}\le\frac{|z_1-z_2|}{\sqrt{2}},\ \ z_1,z_2\in D,
$$
and $P_1^{-1}\circ P_1^{-1}(\C)\ss D$. Thus, {\bf Hypothesis 1} is satisfied and we can use Theorem \ref{T1} and its Corollaries. To parameterize zeros of $f$, we should use the set 
$$
 \S=\{\s:\N\to\{\pm1\},\ \lim_{n\to\iy}\s_n=1\}.
$$ 
Then zeros of $f$ have form \er{all_zeros}
$$
 z(\s)=-\prod_{n=1}^{\infty}\frac{4}{2+2\s_n\sqrt{\frac12+...+\frac{\s_1}2\sqrt{\frac12}}}=-\prod_{n=1}^{\infty}\frac{1}{\frac12+\frac{\s_n}2\sqrt{\frac12+...+\frac{\s_1}2\sqrt{\frac12}}}.
$$
Computations show
$$
 z(1,1,1,...)=-\frac{\pi^2}{8},
\ \ 
 z(-1,1,1,...)=-\frac{9\pi^2}{8},
\ \ 
 z(-1,-1,1,...)=-\frac{25\pi^2}{8},
 \ \ 
 z(1,-1,1,...)=-\frac{49\pi^2}{8}
$$
and so on. This is in full agreement with expected values, since $f(z)=\cos\sqrt{-2z}$. In this case, the formulas for zeros are, in fact, modified Vi\`ete's formulas, see also \cite{S,N}. The order of entire function $f$ is $1/2$. WH-factorization is
$$
 \cos\sqrt{-2z}=\prod_{n=1}^{\iy}\lt(1+\frac{8z}{(2n-1)^2\pi^2}\rt)=\prod_{\s\in\S}\lt(1+z\prod_{n=1}^{\iy}\lt(\frac12+\frac{\s_n}2\sqrt{\frac12+...+\frac{\s_1}2\sqrt{\frac12}}\rt)\rt).
$$

{\bf 2.} Consider the case $P(z)=z^2-1$. SP-equation is $f(az)=f(z)^2-1$. We take $f(0)=b=\frac{\sqrt{5}+1}2$, $f'(0)=1$. Then $a=(z^2-1)'|_{z=b}=2b$. Polynomial \er{Q} is
$$
Q(z)=\frac{z^2-1-b}{z-b}=z+b.
$$
There are two branches of $P^{-1}$:
$$
P_{1}^{-1}(w)=\sqrt{1+w},\ \ P_{-1}^{-1}(w)=-\sqrt{1+w}.
$$
Again, using the arguments from {\bf Example 1}, we can state that  {\bf Hypothesis 1} is satisfied. To parametrize zeros of $f$, we should use the same set as in the previous example 
$$
\S=\{\s:\N\to\{\pm1\},\ \lim_{n\to\iy}\s_n=1\}.
$$ 
Then zeros of $f$ have the form
$$
z(\s)=-b\prod_{n=1}^{\infty}\frac{2b}{b+\s_n\sqrt{1+...+\s_1\sqrt{1}}}.
$$
\begin{figure}
	\begin{minipage}[h]{0.32\linewidth}
		\center{\includegraphics[width=0.99\linewidth]{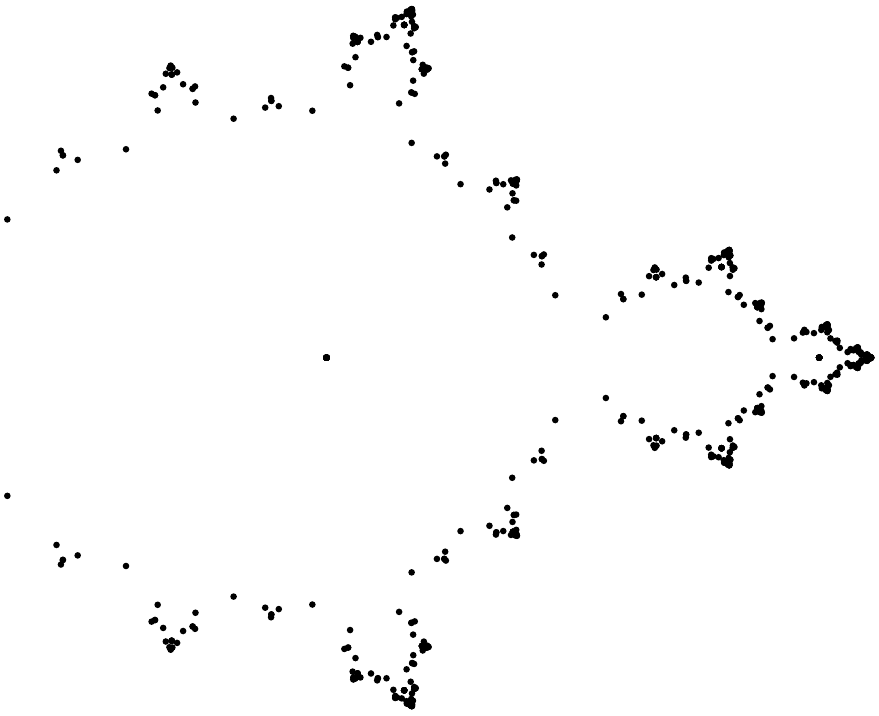} \\ (a)}
	\end{minipage}
	\hfill
	\begin{minipage}[h]{0.32\linewidth}
		\center{\includegraphics[width=0.99\linewidth]{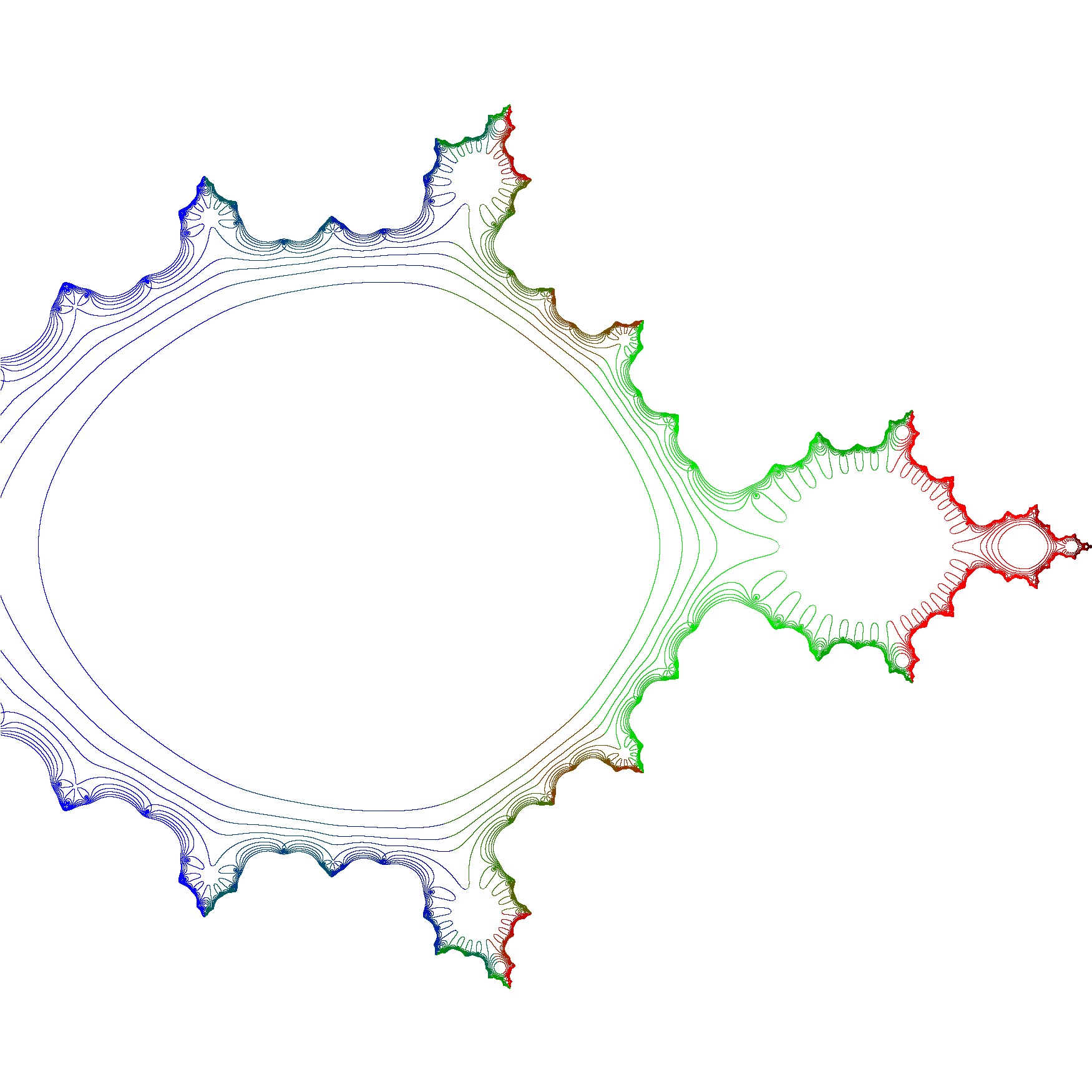} \\ (b)}
	\end{minipage}
    \hfill
    \begin{minipage}[h]{0.32\linewidth}
	    \center{\includegraphics[width=0.99\linewidth]{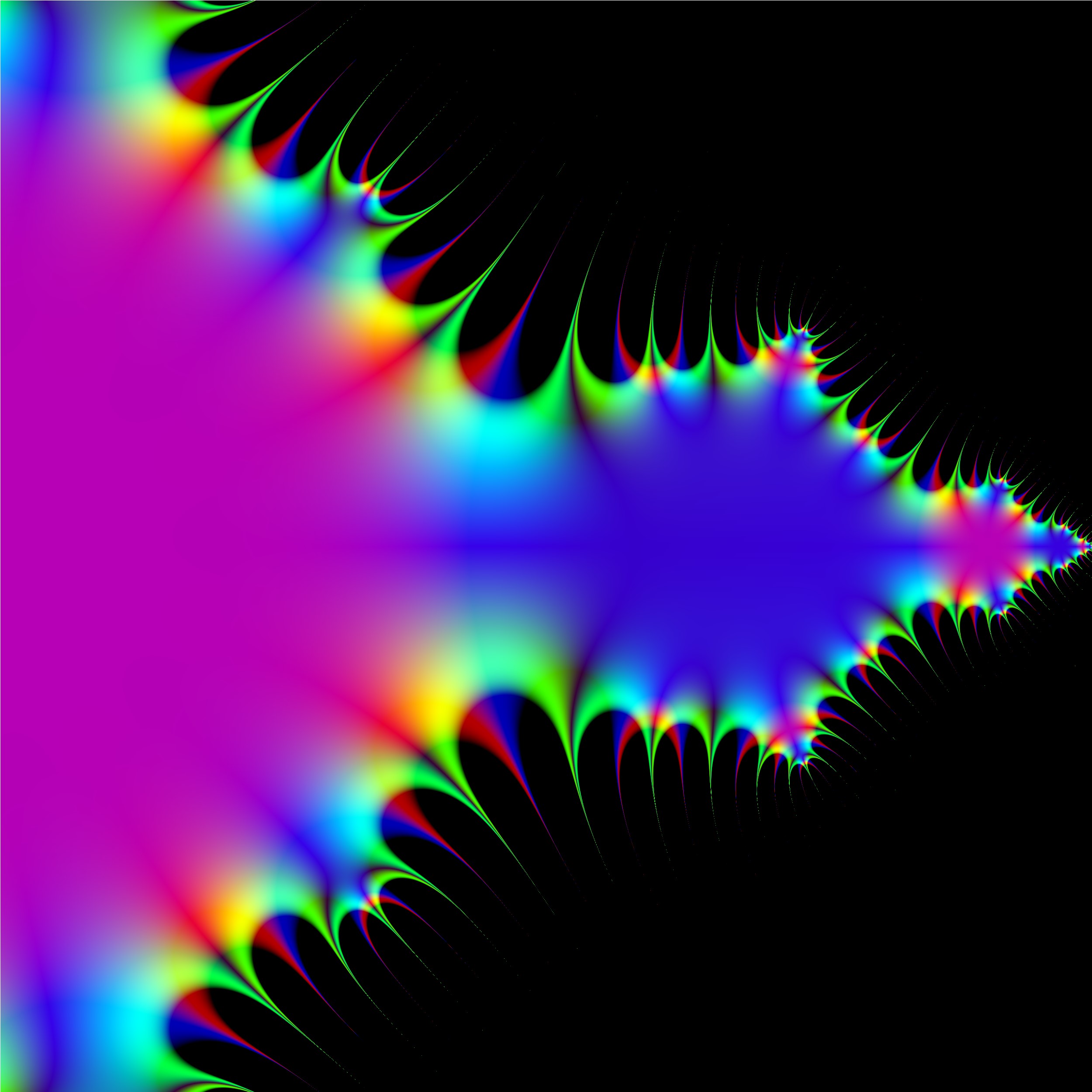} \\ (c)}
    \end{minipage}
	\caption{(a) zeros of $f(z)$ in the complex plane, $|z|\le5\cdot10^5$; (b) images $f^{-1}(w)(=g_{\s}(w))$, see \er{branchesg}, of circles  $|w|=\frac{r^3}{125.05}$, $r=1,...,10$ for the first $10^3$ values $\s\in\S$ depicted in different colors; (c) real and imaginary  parts of $f(z)$, where $\min\{|\Re f(z)|,|\Im f(z)|\}\le2$ and $10^{-3}z\in[-2.5,0]\ts[-1.25,1.25]$.}\lb{fig1}
\end{figure}
The first negative zero $z(1,1,1,...)=-2C$ relates to the so-called Paris constant $C$ appearing in the approximation of the golden ratio by nested square root radicals, see \cite{P,F,K}. Zeros of $f$ are also related to the polynomial dynamics generated by $P=z^2-1$ and, hence, approximate the corresponding Julia set growing up to infinity, see more in \cite{EL,ES,DGV}. The zeros form impressive fractal structures, see Fig. \ref{fig1}. The order of entire function $f$ is $\ln2/\ln a<1$. Hence, there is WH-factorization
$$
f(z)=b\prod_{\s\in\S}\lt(1+\frac{z}{b}\prod_{n=1}^{\infty}\frac{b+\s_n\sqrt{1+...+\s_1\sqrt{1}}}{2b}\rt).
$$
There are infinitely many complex zeros of multiplicities $2^n$ for any $n\ge0$, see \cite{K}. All the multiplicities are taken into account in WH-factorization mentioned above. The first, second and third momentum formulas for zeros, see \er{moment}, \er{moment2} and \er{momentm}, are
$$
\sum_{\s\in\S}\prod_{n=1}^{\infty}\frac{(b+\s_n\sqrt{1+...+\s_1\sqrt{1}})^m}{(2b)^m}=\ca 1,& m=1, \\  1-\frac{1}{\sqrt{5}},& m=2,\\ \frac25,& m=3. \ac
$$

{\bf 3.} Let us consider the cubic SP-equation $f(az)=f(z)^3-6$, $f(0)=b=2$, $f'(0)=1$. Then $a=3b^2=12$. The order of the entire function $f(z)$ is $\ln 3/\ln 12<1$. Let us skip the similar arguments as in the previous examples that show that the principal branch of $P^{-1}$ for $P(z)=z^3-6$ satisfies {\bf Hypothesis 1}. So, we can use \er{moment}, \er{moment2} to obtain explicit momentum formulas
$$
\sum_{k_n\in\{0,1,2\}; \ \lim k_n=0}\ \prod_{n=1}^{\iy}\frac{(e^{\frac{2\pi i k_n}3}\sqrt[3]{6+...+e^{\frac{2\pi i k_1}3}\sqrt[3]{6}})^2+2e^{\frac{2\pi i k_n}3}\sqrt[3]{6+...+e^{\frac{2\pi i k_1}3}\sqrt[3]{6}}+4}{12}=1,
$$
$$
 \sum_{k_n\in\{0,1,2\}; \ \lim k_n=0}\lt(\prod_{n=1}^{\iy}\frac{(e^{\frac{2\pi i k_n}3}\sqrt[3]{6+...+e^{\frac{2\pi i k_1}3}\sqrt[3]{6}})^2+2e^{\frac{2\pi i k_n}3}\sqrt[3]{6+...+e^{\frac{2\pi i k_1}3}\sqrt[3]{6}}+4}{12}\rt)^2=\frac{9}{11}
$$
and so on. All of the momenta are rational numbers.
 
{\section{Proof of Theorems \ref{T1}, \ref{T2}, and \ref{T3}} \lb{Proof}} 

First of all let us show that infinite products \er{all_zeros} are well defined. Suppose that
\[\lb{pf1}
 Q(P_{\s_n}^{-1}\circ...\circ P_{\s_1}^{-1}(w))=0
\]
for some $w\in\C$ and $n\in\N$. If $n=1$ then \er{Q} gives $w=b$. Consider the case $n>1$. We have that $P_{\s_n}^{-1}\circ...\circ P_{\s_1}^{-1}(w)\ne b$, since $Q(b)=P'(b)=a\ne0$. Next, if $P_{\s_n}^{-1}\circ...\circ P_{\s_1}^{-1}(w)\ne b$ then \er{Q} and \er{pf1} give us
$$
 P(P_{\s_n}^{-1}\circ...\circ P_{\s_1}^{-1}(w))=P_{\s_{n-1}}^{-1}\circ...\circ P_{\s_1}^{-1}(w)=b,
$$
which leads to $w=b$, since $P(b)=b$. Hence any denominator in \er{all_zeros} is non-zero, since $b\ne0$ by the assumption from the beginning of the article. Due to analyticity of $P_0^{-1}$ in some open neighbourhood of its attracting point $b$, where $|(P_0^{-1})'(b)|=|a^{-1}|<1$, we have that for any $\ve>0$ there is $\d>0$ such that for any $w\in\{|w-b|<\d\}$:
\[\lb{Remark}
 P_0^{-1}(b+w)=b+R(w),\ \ |R(w)|<(|a^{-1}|+\ve)|w|.
\]
Identity and inequality \er{Remark} along with \er{H1} and the stability condition $\lim\s_n=0$ in \er{S} lead to
\[\lb{e0}
 P_{\s_{n}}^{-1}\circ...\circ P_{\s_1}^{-1}(0)=b+O(c^{n}),\ \ n\to\iy,
\] 
where small $\ve>0$ is taken such that $c:=|a^{-1}|+\ve<1$. Hence
$$
 Q(P_{\s_{n}}^{-1}\circ...\circ P_{\s_1}^{-1}(0))=Q(b)+O(c^{n})=P'(b)+O(c^{n})=a+O(c^{n}).
$$ 
This guaranties the convergence of infinite products \er{all_zeros} for any $\s\in\S$. The convergence of the products is exponentially fast, since $c<1$. Note that the same arguments give also the exponential rate of convergence of \er{g0p}.

Let $\wt z$ be some zero of $f(z)$ of multiplicity $\wt m\in\N$, i.e. 
\[\lb{df}
 f^{(j)}(\wt z)=0\ \ {\rm for}\ \  j=0,...,\wt m-1. 
\]
SP-equation gives
\[\lb{e1}
 f(z)=P^{\circ n}(f(a^{-n}z)),\ \ n\in\N.
\]
Taking $n>\wt m$ such that $f'(a^{-n}\wt z)\ne0$ (recall that $f'(0)=1\ne0$ and $|a|>1$), differentiating \er{e1} at $z=\wt z$ and using \er{df}, we obtain that
\[\lb{e2}
 (P^{\circ n})^{(j)}(f(a^{-n}\wt z))=0,\ \ j=0,...,\wt m-1.
\]
Thus, $f(a^{-n}\wt z)$ is a root of $P^{\circ n}$ of a multiplicity at least $\wt m$. This means that
\[\lb{e3}
 f(a^{-n}\wt z)=P_{\s_n}^{-1}\circ...\circ P_{\s_1}^{-1}(0)
\]
for at least $\wt m$ different $\s:\{1,...,n\}\to\{0,...,d\}$. 

SP-equation can be written in the form $f(a^{-1}z)=P^{-1}(f(z))$. Since $f(0)=b$, the branch $P^{-1}$ should coincide with the principal branch $P_0^{-1}$ in a small neighbourhood of $b$, i.e. $f(a^{-1}z)=P_0^{-1}(f(z))$ for all sufficiently small $z$ (see also the remark before {\bf Hypothesis 1}). Let $\wt n$ be such that $f(a^{-\wt n}\wt z)$ belongs to this small neighbourhood of $b$. We assume also that $\wt n$ is large enough to satisfy \er{e3} with at least $\wt m$ different $\s$. Then, by \er{e3}, we have
\[\lb{e4}
 f(a^{-\wt n-k}\wt z)=\underbrace{P_0^{-1}\circ...\circ P_0^{-1}}_k\circ P_{\s_{\wt n}}^{-1}\circ...\circ P_{\s_1}^{-1}(0),\ \ k\ge0.
\]
Denote $\s_n=0$ for $n>\wt n$. Thus, using $f(0)=b$, $f'(0)=1$, we get
\begin{multline}
 \wt z=\lim_{n\to\iy}a^n(f(a^{-n}\wt z)-b)=\lim_{n\to\iy}a^n(P_{\s_n}^{-1}\circ...\circ P_{\s_1}^{-1}(0)-b)=\\
 \lim_{n\to\iy}\frac{a(P_{\s_n}^{-1}\circ...\circ P_{\s_1}^{-1}(0)-b)}{P_{\s_{n-1}}^{-1}\circ...\circ P_{\s_1}^{-1}(0)-b}a^{n-1}(P_{\s_{n-1}}^{-1}\circ...\circ P_{\s_1}^{-1}(0)-b)=\\
 \lim_{n\to\iy}\frac{a}{Q(P_{\s_n}^{-1}\circ...\circ P_{\s_1}^{-1}(0))}(P_{\s_{n-1}}^{-1}\circ...\circ P_{\s_1}^{-1}(0)-b)=-b\prod_{n=1}^{\iy}\frac{a}{Q(P_{\s_n}^{-1}\circ...\circ P_{\s_1}^{-1}(0))}.\lb{e5}
\end{multline}
Like \er{e3}, identity
\[\lb{e6}
 \wt z=-b\prod_{n=1}^{\iy}\frac{a}{Q(P_{\s_n}^{-1}\circ...\circ P_{\s_1}^{-1}(0))}
\]
holds for at least $\wt m$ different $\s\in\S$. 

Conversely, suppose that \er{e6} holds for $\wt m$ different $\s\in\S$. To finish the proof we need to show that $\wt z$ is a zero of $f$ of a multiplicity at least $\wt m$. Using \er{limit} and the second identity in \er{e5}, we obtain
\begin{multline}
 f(\wt z)=\lim_{n\to\iy}f(a^n(P_{\s_n}^{-1}\circ...\circ P_{\s_1}^{-1}(0)-b))=\\
 \lim_{n\to\iy}\underbrace{P\circ...\circ P}_n(b+a^{-n}(a^n(P_{\s_n}^{-1}\circ...\circ P_{\s_1}^{-1}(0)-b)))=\\
 \lim_{n\to\iy}\underbrace{P\circ...\circ P}_n\circ P_{\s_n}^{-1}\circ...\circ P_{\s_1}^{-1}(0)=0,\lb{z1}
\end{multline}
since the convergence of \er{limit} is uniform in any bounded domain. Hence $\wt z$ is a zero of $f$. Now, let $N$ be such that $P_{\s_N}^{-1}\circ...\circ P_{\s_1}^{-1}(0)$ is sufficiently close to $b$, where $f^{-1}$ is defined, so that
\[\lb{f1}
 P_{\s_N}^{-1}\circ...\circ P_{\s_1}^{-1}(0)=f(t_{\s}).
\]
We can do this because $P_{\s_n}^{-1}=P_0^{-1}$ for large $n$ by definition \er{S}, and we are under {\bf Hypothesis 1}. We also assume that $N$ is so large that \er{f1} is valid for at least $\wt m$ different $\s:\{1,...,N\}\to\{0,...,d\}$ coinciding with the segments of those $\s\in\S$ mentioned in \er{e6}, and, also, all $\s_n=0$ for $n>N$. 
For simplicity, in the previous sentence we use the same symbol for $\s\in\S$ and for its segment $\s|_{\{1,...,N\}}$. 
Finally, it is assumed that $N$ is so large that
\[\lb{f2}
 \underbrace{P_0^{-1}\circ...\circ P_0^{-1}}_n\circ P_{\s_N}^{-1}\circ...\circ P_{\s_1}(0)=f(a^{-n}t_{\s}),
\]  
see comments before \er{e4}. Using \er{f2}, \er{e6}, the assumption $\s_n=0$, $n>N$, and the arguments similar to \er{e5}, we obtain
\[\lb{f3}
 t_{\s}=\lim_{n\to\iy}a^n(f(a^{-n}t_{\s})-b)=-a^{-N}b\prod_{n=1}^{\iy}\frac{a}{Q(P_{\s_n}^{-1}\circ...\circ P_{\s_1}^{-1}(0))}=a^{-N}\wt z.
\]
Hence, all $t_{\s}$ are equal to each other. Using \er{f1}, the remark after \er{f1} about $\wt m$ different $\s$, and \er{f3}, we conclude that $f(a^{-N}\wt z)$ is a zero of $P^{\circ N}$ of a multiplicity at least $\wt m$. Thus, differentiating $f(z)=P^{\circ N}(f(a^{-N}z))$ at $z=\wt z$, we obtain that $f^{(j)}(\wt z)=0$, $j=0,...,\wt m-1$. Hence, $\wt z$ is a zero of $f$ of a multiplicity at least $\wt m$. The proof of Theorem \ref{T1} is finished.

To prove Theorem \ref{T2} let us note that
\[\lb{gtail}
 g_{\s}(w)=a^m(P_{\s_m}^{-1}\circ...\circ P_{\s_1}^{-1}(w)-b)\prod_{n=m+1}^{\iy}\frac{a}{Q(P_{\s_n}^{-1}\circ...\circ P_{\s_1}^{-1}(w))},\ \ m\ge0,\ \ w\ne b,
\]
where $P_{\s_0}^{-1}(w):=w$ and $g_{\s}(w)$ are defined in \er{branchesg}. Moreover, using the similar arguments as in the proof of Theorem \ref{T1}, we can state that if $\s_m\ne0$ then \er{gtail} is true for $w=b$. Now, choosing the first $m$ such that $\s_m\ne0$, we finish the proof of Theorem \ref{T2}.

Recall that $P_0^{-1}$ is the unique branch of $P^{-1}$ such that $P_0^{-1}(b)=b$ and $P_0^{-1}$ is analytic in the neighbourhood of $b$. Let $\d>0$ be so small that 
\[\lb{bcirc}
 P_{j}^{-1}(\C)\cap\{w:\ |w-b|<\d\}=\es,\ \ \forall j\ne0.
\]
Identities \er{branchesg} and \er{gtail} give
\[\lb{gc}
 f(a^{-m}g_{\s}(w))=P_{\s_m}^{-1}\circ...\circ P_{\s_1}^{-1}(w).
\]
If $\s_m\ne0$ then \er{bcirc} leads to 
\[\lb{gc1}
 |f(a^{-m}g_{\s}(w))-b|\ge\d.
\]
Since $f$ is continuous and $f(0)=b$, there is $C>0$ such that 
\[\lb{gc2}
 \{z:\ |z|<C\}\ss f^{-1}(\{w:\ |w-b|<\d\}).
\]
Thus, by \er{gc1} and \er{gc2}, we have $|a^{-m}g_{\s}(w)|\ge C$, which gives \er{est1}. Using \er{est1}, we get
\[\lb{gc3}
 \sum_{|\supp\s|\ge N}|g_{\s}(w)|^{-m}=\sum_{n\ge N}\sum_{|\supp\s|=n}|g_{\s}(w)|^{-m}\le\frac1{C^m}\sum_{n\ge N}d^n|a|^{-mn}=O((d|a|^{-m})^N).
\]
By analogy, we can estimate the product in \er{est2}. The proof of Theorem \ref{T3} is finished.

\end{document}